# STABILIZATION OF THE WITT GROUP

## by Max KAROUBI

**Abstract.**

In this Note, using an idea due to Thomason [8], we define a "homology theory" on the category of rings which satisfies excision, exactness, homotopy (in the algebraic sense) and periodicity of order 4. For regular noetherian rings, we find Balmers's higher Witt groups. For more general rings, this homology is isomorphic to the KT-theory of Hornbostel [3], inspired by the work of Williams [9]. For real or complex C*-algebras, we recover - up to 2 torsion - topological K-theory.

**1.** Let A be a ring with an antiinvolution $a \mapsto \bar{a}$ and let $\varepsilon$ be an element of the center of A such that $\varepsilon\bar{\varepsilon} = 1$. We assume also that 2 is invertible in the ring. There are now well known definitions of the higher hermitian K-group (denoted by $_{\varepsilon}L_n(A)$, as in [5]) and the higher Witt group $_{\varepsilon}W_n(A)$ : this is the cokernel of the map induced by the hyperbolic functor

$$K_n(A) \longrightarrow \, _{\varepsilon}L_n(A)$$

where the $K_n(A)$ denote the Quillen K-group (which is defined for all values of $n \in \mathbf{Z}$).

One of the fundamental results of higher Witt theory is the periodicity isomorphism (where $\mathbf{Z}' = \mathbf{Z}[1/2]$, cf.[4])

$$_{\varepsilon}W_n(A) \otimes \mathbf{Z}' \cong \, _{-\varepsilon}W_{n-2}(A) \otimes \mathbf{Z}'$$

It is induced by the cup-product with a genuine element $u_2 \in \, _{-1}L_{-2}(\mathbf{Z}')$. By analogy with algebraic topology, we shall call $u_2$ the Bott element in Witt theory. This element is explicitly described in the following way. We consider the 2 x 2 matrix (with the involution defined by $\bar{z} = z^{-1}$ and $\bar{t} = t^{-1}$ and where we put $\lambda = \bar{\lambda} = 1/2$).

$$M = \begin{pmatrix} b(t+t^{-1}-2) & a + (1-a)t \\ -(1-d+dt^{-1}) & -c \end{pmatrix} \text{ where } \begin{vmatrix} a & b \\ c & d \end{vmatrix} = u\, p_0\, u^{-1}$$

$$\text{with } p_0 = \begin{vmatrix} 1 & 0 \\ 0 & 0 \end{vmatrix} \text{ and } u = \begin{pmatrix} \bar{\lambda}z + \lambda & \lambda\bar{\lambda}(z-1) \\ z-1 & \lambda z + \bar{\lambda} \end{pmatrix}$$

This 2 x 2 matrix represents an element of $_{-1}L_0(\mathbf{Z}'[t, t^{-1}, z, z^{-1}])$ whose image in $_{-1}L_{-2}(\mathbf{Z}') \cong \mathbf{Z} \oplus \mathbf{Z}/2$ is a free generator (cf. [5] for the details).

**2**. The higher Witt groups $_\varepsilon W_n(A)$ do not have all the nice formal properties one should expect. For instance, a cartesian square of rings with antiinvolutions (where the vertical maps are surjective)

$$\begin{array}{ccc} A & \longrightarrow & A_1 \\ \downarrow & & \downarrow \\ A_2 & \longrightarrow & A' \end{array}$$

does not induce in general a long Mayer-Vietoris exact sequence of Witt groups

$$\longrightarrow {}_\varepsilon W_{n+1}(A') \longrightarrow {}_\varepsilon W_n(A) \longrightarrow {}_\varepsilon W_n(A_1) \oplus {}_\varepsilon W_n(A_2) \longrightarrow {}_\varepsilon W_n(A') \longrightarrow$$

As a counterexample for n = 0 coming from topology, one might take the ring of complex continuous functions on finite CW-complexes, provided with the trivial involution (otherwise this will imply - with topological notations - that the obvious cokernel from KU to KO is a cohomology theory and therefore a direct factor in KO). Note however that by tensoring with **Z'**, we restore the Mayer-Vietoris axiom as a consequence of the periodicity theorem and the direct splitting of hermitian K-theory shown in [4] p. 253.

Following an idea due to Thomason [8], one may overcome this difficulty by **stabilizing** the higher Witt groups. More precisely, we define a new theory $_\varepsilon \mathcal{W}_n(A)$ as the limit of the inductive system

$$_\varepsilon W_n(A) \longrightarrow {}_{-\varepsilon} W_{n-2}(A) \longrightarrow {}_\varepsilon W_{n-4}(A) \longrightarrow \ldots$$

where the arrows are induced by the cup-product with the Bott element $u_2$ mentioned above. As a matter of fact, the periodicity map $_\varepsilon W_n(A) \longrightarrow {}_{-\varepsilon} W_{n-2}(A)$ can be factored as

$$_\varepsilon W_n(A) \longrightarrow {}_{-\varepsilon} L_{n-2}(A) \longrightarrow {}_{-\varepsilon} W_{n-2}(A)$$

Therefore $_\varepsilon \mathcal{W}_n(A)$ is also the limit of the inductive system

$$_\varepsilon L_n(A) \longrightarrow {}_{-\varepsilon} L_{n-2}(A) \longrightarrow {}_\varepsilon L_{n-4}(A) \longrightarrow \ldots$$

**3. THEOREM.** *This new theory $_\varepsilon \mathcal{W}_n(A)$ satisfies the following properties*
**a) Homotopy invariance**. *The polynomial extension* $A \longrightarrow A[t]$, *where* $\bar{t} = t$, *induces an isomorphism*

$$_\varepsilon \mathcal{W}_n(A) \cong {}_\varepsilon \mathcal{W}_n(A[t])$$

**b) Exactness and excision**. *From a cartesian square as above (with $\psi$ surjective)*

$$\begin{array}{ccc} A & \longrightarrow & A_1 \\ \phi \downarrow & & \psi \downarrow \\ A_2 & \longrightarrow & A' \end{array}$$

one deduces an isomorphism of the associated relative groups

$$_\varepsilon \mathcal{W}_n(\phi) \cong {}_\varepsilon \mathcal{W}_n(\psi)$$

and therefore a Mayer-Vietoris *exact sequence* (*for all* $n \in \mathbb{Z}$)

$$\longrightarrow {}_\varepsilon \mathcal{W}_{n+1}(A') \longrightarrow {}_\varepsilon \mathcal{W}_n(A) \longrightarrow {}_\varepsilon \mathcal{W}_n(A_1) \oplus {}_\varepsilon \mathcal{W}_n(A_2) \longrightarrow {}_\varepsilon \mathcal{W}_n(A') \longrightarrow$$

**c) Periodicity**. *The cup-product with the* Bott *element induces the isomorphisms*

$$_\varepsilon \mathcal{W}_n(A) \cong {}_{-\varepsilon} \mathcal{W}_{n-2}(A) \cong {}_\varepsilon \mathcal{W}_{n-4}(A)$$

**d) Normalization.** *Let us assume now that* A *is a regular noetherian ring. Then the the natural map*

$$_\varepsilon W_0(A) \longrightarrow {}_\varepsilon \mathcal{W}_0(A)$$

*is an isomorphism and the group* $_\varepsilon \mathcal{W}_1(A)$ *is isomorphic to the cokernel of the map defined in* [5]

$$k_0(A) \longrightarrow {}_\varepsilon W_1(A)$$

*Moreover, the groups* $_1\mathcal{W}_n(A)$ *coincide with the higher* Witt *groups of* Balmer [1].

*Proof.* Periodicity is imposed by the definition (as in Thomason's theory). Homotopy invariance is a consequence of the same property for the Witt groups. Since the $L_n$-groups satisfy the excision and exactness properties for $n < 0$ (cf. [6] for instance), this is also true of the theory $_\varepsilon \mathcal{W}_n$ : as we have noticed before, $_\varepsilon \mathcal{W}_n(A)$ is also the limit of the inductive system

$$_\varepsilon L_n(A) \longrightarrow {}_{-\varepsilon} L_{n-2}(A) \longrightarrow {}_\varepsilon L_{n-4}(A) \longrightarrow$$

If A is regular noetherian, the K-theory groups $K_n(A)$ are 0 for $n < 0$. Therefore, according to the 12 term exact sequence proved in [5], we have an exact sequence

$$0 = k_{-1}(A) \longrightarrow {}_\varepsilon W_0(A) \longrightarrow {}_{-\varepsilon} W'_{-2}(A) \longrightarrow k'_{-1}(A) = 0$$

with an obvious isomorphism $_{-\varepsilon}W'_{-2}(A) \cong {}_{-\varepsilon}W_{-2}(A)$ since again $K_n(A) = 0$ for $n = -1$ and $-2$. With the same argument, we prove that $_\varepsilon W_n(A) \cong {}_{-\varepsilon}W_{n-2}(A)$ for all $n \le 0$.
In the same spirit, we have an exact sequence

$$k_0(A) \longrightarrow {}_\varepsilon W_1(A) \longrightarrow {}_\varepsilon \mathcal{W}_1(A) \longrightarrow 0$$

The first map is the following. An element of $k_0(A)$ is the class of a module E which is isomorphic to its dual. Its image in $_\varepsilon W_1(A)$ is associated to the automorphism of the hyperbolic module $E \oplus E^* \cong E \oplus E$ defined by the matrix

$$\begin{pmatrix} 0 & 1 \\ \varepsilon & 0 \end{pmatrix}$$

Finally, the isomorphism with the higher Witt groups defined by Balmer will follow from next theorem.

**4. THEOREM.** *The homology $_\varepsilon\mathbb{W}_*(A)$ is isomorphic to the $KT_*$-theory of Hornbostel (cf. [3] § 5 ).*

*Proof.* This KT-theory is the direct limit of the system

$$_\varepsilon L_n(A) \longrightarrow {}_\varepsilon L_n(U_A) \longrightarrow \ldots \longrightarrow {}_\varepsilon L_n(U_A^r) \longrightarrow$$

where $U_A$ is the ring defined in [5], p. 263 and $U_A^r$ the r-iteration of the "U-construction". All the arrows above are $L_*$-module maps as defined in [4] p. 233 and [5] p. 276. This implies that the homomorphism

$$_\varepsilon L_n(A) \longrightarrow {}_\varepsilon L_n(U_A^r)$$

is the cup-product with a well defined element $w_r$ in $_1L_0(U_{\mathbf{Z}'}^r)$ (this is probably related to the question 6.6 raised by Hornbostel in his paper [3]). On the other hand, as a consequence of the fundamental theorem of hermitian K-theory (cf. [5] p. 264), we have an isomorphism of $L_*$-modules between $_\varepsilon L_n(U_A^r)$ and $_\varepsilon L_n(U_A^{r+4})$ (as noticed also by Williams [9]). Therefore, the previous direct limit is simply the limit of the system

$$_\varepsilon L_n(A) \longrightarrow {}_\varepsilon L_{n-4}(A) \longrightarrow \ldots$$

where the arrows are defined by the cup-product with a specific element w in $_1L_{-4}(\mathbf{Z}')$. On the other hand, we know that if we apply this construction to the ring $A = \mathbf{Z}'$ and $\varepsilon = 1$, we find an isomorphism between $_1W_0(\mathbf{Z}')$ and $_1L_{-4}(\mathbf{Z}') = {}_1W_{-4}(\mathbf{Z}')$ (because the ring $\mathbf{Z}'$ is regular). As a matter of fact, we find a chain of isomorphisms

$$_1W_0(\mathbf{Z}') \cong {}_1L_0(U_{\mathbf{Z}'}) \cong {}_1L_0(U_{\mathbf{Z}'}^2) \cong {}_1L_0(U_{\mathbf{Z}'}^3) \cong {}_1L_0(U_{\mathbf{Z}'}^4) \cong {}_1L_{-4}(\mathbf{Z}')$$

I claim that w, the image of 1 by this chain of isomorphisms, is $(u_2)^2$ up to a unipotent element. This is exactly the well known computation of the classical Witt ring of $\mathbf{Z}'$ which is $\mathbf{Z} \oplus \mathbf{Z}/2$ generated by the classes of the following elements in the Grothendieck Witt group : $<x^2>$ and $<x^2> - <2x^2>$.

If A is regular noetherian and $\varepsilon = 1$, Hornbostel has proved moreover in [3] that $KT_n(A)$ is isomorphic to the n-Witt group defined by Balmer, which proves the last part of theorem 3

**5. Remark.** For simplicity's sake, we have just considered hermitian K-theory <u>groups</u>. One could have taken as well homotopy colimits of the corresponding classifying spaces, using for instance the machinery developped in [4] § 1, greatly generalized by Schlichting.

**6. THEOREM.** *The theory $_\varepsilon W_n(A)$ is invariant under nilpotent extensions. In other words, let* I *be a nilpotent ideal of* A, *stable under the antiinvolution. Then the quotient map* A $\longrightarrow$ A/I *induces an isomorphism*

$$_\varepsilon W_n(A) \cong {_\varepsilon W_n(A/I)}$$

*Proof.* For $n < 0$, $_\varepsilon L_n(A)$ is $_\varepsilon L_0$ of the $(-n)^{th}$ suspension $S^{-n}A$. Therefore, it suffices to show that the induced map $_\varepsilon L_0(A) \longrightarrow {_\varepsilon L_0(A/I)}$ is an isomorphism. To prove surjectivity, we remark that every A/I-hermitian module M is the image of a self-adjoint projection operator Q in some hyperbolic module $H(B^n)$, where $B = A/I$ (cf. [6]). We write $Q = (J - 1)/2$ where J is an involution. We lift J to an operator on $H(A^n)$ which we call R. Then $R^* = R + \eta$, where $\eta$ is a matrix with coefficients in I. By replacing R by $S = R + \eta/2$, we see that we may assume R to be self-adjoint. Now $R^2 = 1 + \gamma$ where $\gamma$ is a matrix with coefficients in I and the power expansion $U = (1 + \gamma)^{-1/2}$ is convergent since I is nilpotent. Therefore, the product RU is a self-adjoint involution which is a lift of J. This proves the surjectivity of the map $_\varepsilon L_0(A) \longrightarrow {_\varepsilon L_0(A/I)}$ (take the image of 1 - RU).

For the proof of injectivity, the argument is quite similar. Let $E_1$ et $E_2$ be two hermitian modules over the ring A which become isomorphic over the ring A/I. This means (up to stabilization) that $E_1$ and $E_2$ are associated[1] to self-adjoint involutions $J_1$ and $J_2$ which are conjugate mod I. In other words, we have $J_2 = \alpha J_1 (\alpha)^{-1}$ where $\alpha$ is a unitary matrix mod. I. This invertible matrix can be lifted to an invertible matrix $\beta$ in $H(A^n)$ and we put

$$\gamma = \beta (\beta\beta^*)^{-1/2}$$

(polar decomposition of matrices). This is a unitary matrix which is a lift of $\alpha$. Therefore, by replacing $J_1$ by $\gamma J_1 \gamma^1$, we may assume that $J_1 = J_2$ mod.I

We now consider the matrix $\delta = (1 + J_1 J_2)/2$. Since $\delta = 1$ mod. I, $\delta$ is an invertible matrix such that $\delta.J_1 = J_2.\delta$. We consider again the polar decomposition of $\delta$, i.e. we replace $\delta$ by $\delta' = \delta (\delta.\delta^*)^{-1/2}$. Since $\delta.\delta^*$ commutes with $J_1$, we also have $\delta'.J_1 = J_2.\delta'$. This shows that $\delta'$ is a unitary matrix which conjugates $J_1$ and $J_2$ and completes the proof of injectivity.

**7. Remark.** Using [2] and most of the above properties of the theory $_\varepsilon W_n$, Schlichting was able to prove cdh descent for the theory $_\varepsilon W_*$, extended to the category of (commutative) schemes of finite type over a field of characteristic 0 [7].

**8.** When A is a Banach algebra, we may consider the topological analogs of the previous definitions. In that case, the group $_\varepsilon \mathcal{W}_n^{top}(A)$ is simply isomorphic to $_\varepsilon W_n^{top}(A) \otimes \mathbf{Z}'$. We prove this fact by looking at the image of $(u_2)^4$ in $_1W_{-8}^{top}(\mathbf{R}) \cong \mathbf{Z}$. According to [5] we find 8 times the generator. Therefore, by taking the inductive limit we localize with respect to the multiplicative system $(2^r)$ ; the natural map $_\varepsilon W_n^{top}(A) \longrightarrow {_\varepsilon \mathcal{W}_n^{top}(A)}$ coincides with this localization. One should also notice that $_\varepsilon W_n^{top}(A) \otimes \mathbf{Z}'$ is isomorphic to $_\varepsilon W_n(A) \otimes \mathbf{Z}'$, as a consequence of the periodicity theorem and the well-known computation of $_\varepsilon W_0$ and $_\varepsilon W_1$. Finally, if A is a C*algebra, it is well known that $_1W_n^{top}(A)$ is isomorphic to the topological K-theory of A.

---

[1] More precisely, $E_1$ (resp. $E_2$) is the image of $(J_1 - 1)/2$ (resp. $(J_2 - 1)/2$).

**9. Aknowledgements**. It is a pleasure to thank Jens Hornbostel et Marco Schlichting for their interesting remarks and suggestions after a first draft of this Note.